\documentclass[11pt,a4paper]{amsart}
\usepackage[T1]{fontenc}
\usepackage{lmodern}
\usepackage[margin=28mm,headheight=14pt]{geometry}
\usepackage{microtype,mathtools,amssymb,amscd,mathrsfs,enumitem,booktabs}
\usepackage[dvipsnames]{xcolor}
\usepackage{hyperref}
\hypersetup{colorlinks=true,linkcolor=MidnightBlue,citecolor=MidnightBlue,urlcolor=MidnightBlue,
 pdftitle={Quasi-F-split heights and derived Morita invariance},
 pdfauthor={Keiho Matsumoto}}
\allowdisplaybreaks[1]
\setlist[enumerate]{leftmargin=2.1em,label=\textup{(\roman*)}}
\setlist[itemize]{leftmargin=1.8em}
\newtheorem{theorem}{Theorem}[section]
\newtheorem{proposition}[theorem]{Proposition}
\newtheorem{lemma}[theorem]{Lemma}
\newtheorem{corollary}[theorem]{Corollary}
\newtheorem{maintheorem}{Theorem}

\theoremstyle{definition}
\newtheorem{definition}[theorem]{Definition}
\newtheorem{example}[theorem]{Example}
\newtheorem{remark}[theorem]{Remark}
\newtheorem{question}[theorem]{Question}

\newcommand{\kk}{k}
\newcommand{\Fp}{\mathbf F_p}

\newcommand{\OO}{\mathcal O}
\newcommand{\TT}{\mathcal T}
\newcommand{\UU}{\mathcal U}
\newcommand{\CC}{\mathcal C}
\newcommand{\XX}{\mathscr X}
\newcommand{\cQ}{\mathcal Q}
\newcommand{\Perf}{\operatorname{Perf}}
\newcommand{\THH}{\operatorname{THH}}
\newcommand{\TR}{\operatorname{TR}}
\newcommand{\HH}{\operatorname{HH}}
\newcommand{\Hom}{\operatorname{Hom}}
\newcommand{\RHom}{\operatorname{RHom}}
\newcommand{\Map}{\operatorname{Map}}
\newcommand{\Spec}{\operatorname{Spec}}
\newcommand{\id}{\operatorname{id}}
\newcommand{\res}{\operatorname{res}}
\newcommand{\cofib}{\operatorname{cofib}}
\newcommand{\colim}{\operatorname*{colim}}
\newcommand{\tr}{\operatorname{tr}}
\newcommand{\R}{\mathbf R}
\newcommand{\LL}{\mathbf L}
\newcommand{\hgt}{\operatorname{ht}}
\newcommand{\Sdim}{\operatorname{Sdim}}
\newcommand{\hm}{\mathfrak m}
\newcommand{\stacks}[1]{\href{https://stacks.math.columbia.edu/tag/#1}{Tag #1}}
\numberwithin{equation}{section}
\title[Quasi-\(F\)-split heights and derived Morita invariance]
{Quasi-\(F\)-split heights and derived Morita invariance}
\author{Keiho Matsumoto}

\begin{document}
\begin{abstract}
In this paper, we study quasi-\(F\)-splitting from the viewpoint of non-commutative algebraic geometry. Let \(k\) be a perfect field of characteristic \(p>0\). For a proper dg-category over \(k\) with a Serre functor, we define quasi-\(F^e\)-split heights for \(e\geq1\) by using finite topological restriction homology with coefficients in a shift of the Serre functor. Using the functoriality of polygonic topological Hochschild homology for laced categories due to Harpaz--Nikolaus--Saunier, we show that these heights are invariant under derived Morita equivalence. For a proper Gorenstein scheme of pure dimension over \(k\), the height for \(e=1\) of its dg-category of perfect complexes coincides with the quasi-\(F\)-split height in the sense of Yobuko. For a complete Noetherian Gorenstein local \(k\)-algebra with residue field \(k\), the dg-category of perfect complexes supported at the closed point recovers the \(n\)-quasi-\(F^e\)-splitting of Tanaka--Witaszek--Yobuko for all \(n,e\geq1\). The proofs use the computation of the zeroth homotopy groups of TR with coefficients due to Dotto--Krause--Nikolaus--Patchkoria, Zariski descent, and Serre duality or local duality.
\end{abstract}
\maketitle
\tableofcontents
\clearpage

\section{Introduction}

The aim of this paper is to construct an analogue of the quasi-\(F\)-split height for non-commutative algebraic varieties, and to compare it with the classical invariants of proper Gorenstein schemes and complete Gorenstein local rings.

Let \(\kk\) be a perfect field of characteristic \(p>0\). Quasi-\(F\)-splitting was introduced by Yobuko \cite{Yob19} as a generalization of \(F\)-splitting, in which a splitting of the Frobenius map is replaced by a splitting involving Witt vectors of finite length. It has applications in birational geometry \cite{KTTWYY25}, and the variant with higher Frobenius exponents, quasi-\(F^e\)-splitting, was introduced by Tanaka--Witaszek--Yobuko \cite{TWY24}. Fedder-type criteria make these notions computable for hypersurfaces \cite{KTY22,Yos25}.

The definition of quasi-\(F\)-splitting involves the structure sheaf, the module structure on Witt vectors and a specific map \(\Phi_{X,n}\) (see Section~\ref{sec:classical}), whereas a general non-commutative algebra has no ring endomorphism \(a\mapsto a^p\).

In this paper, we use finite topological restriction homology \(\TR^n\), which has Frobenius and restriction maps, in place of Witt vectors. In order to obtain the splitting condition for a proper scheme from its dg-category of perfect complexes, we take coefficients in the Serre functor. For a Gorenstein scheme, this coefficient corresponds to the canonical bundle. Then the lowest homotopy group of a certain pushout of TR spectra is identified with the top cohomology of the classical pushout \(Q_{X,n}\) tensored with the canonical bundle, and Serre duality translates the injectivity of the induced map into the existence of the required splitting.

\subsection{The categorical construction}
Throughout this paper, all dg-categories are assumed to be small, idempotent-complete and \(\kk\)-linear, and to admit a classical generator. Equivalently, such a dg-category is of the form \(\Perf(A)\) for a dg-algebra \(A\) over \(\kk\). Recall that a dg-category \(\TT\) is proper if the total cohomology of \(\Hom_\TT(x,y)\) is finite-dimensional for all objects \(x,y\). Let \(\TT\simeq\Perf(A)\) be a proper dg-category. We say that \(\TT\) admits a Serre functor if the dual bimodule
\[
 A^\vee:=\RHom_{\kk}(A,\kk)
\]
is invertible with respect to the derived tensor product. In this case, \(S_\TT:=-\otimes^{\LL}_AA^\vee\) is a Serre functor of \(\TT\), that is, there are natural equivalences \(\RHom_\TT(x,y)^\vee\simeq\RHom_\TT(y,S_\TT x)\), and we always regard \(S_\TT\) as a dg-functor in this way. These notions do not depend on the choice of \(A\).

Fix an integer \(d\), and put \(P_{\TT}:=S_{\TT}[-d]\), which we regard as a coefficient bimodule. We use the polygonic spectra of Krause--McCandless--Nikolaus \cite{KMN23}, which refine topological Hochschild homology with coefficients. For an exact endofunctor \(P\) of \(\TT\), we denote by \(\TR^m(\TT;P)\) the finite TR with respect to the truncation set \(\{1,p,\ldots,p^{m-1}\}\); as usual, we have \(\TR^1(\TT;P)=\THH(\TT;P)\). There are restriction and Frobenius maps
\[
 \res\colon\TR^{m+1}(\TT;P)\longrightarrow\TR^m(\TT;P),
 \qquad
 F\colon\TR^{m+1}(\TT;P)\longrightarrow\TR^m(\TT;P^{\otimes p}),
\]
where the Frobenius changes the coefficient from \(P\) to \(P^{\otimes p}\).

For \(n,e\geq1\), we define
\begin{equation}\label{eq:intro-Q}
\cQ^e_n(\TT;P):=
\cofib\left[
\TR^{n+e}(\TT;P)
\xrightarrow{(F^e,\,-\res^{n+e-1})}
\TR^n(\TT;P^{\otimes p^e})\oplus\THH(\TT;P)
\right],
\end{equation}
and we denote by \(\eta^e_n(\TT;P)\colon\THH(\TT;P)\to\cQ^e_n(\TT;P)\) the induced map. We define
\[
\hgt^e_d(\TT):=
\inf\{n\geq1\mid\pi_{-d}\bigl(\eta^e_n(\TT;P_\TT)\bigr)\text{ is injective}\},
\]
where we put \(\inf\emptyset=\infty\). The integer \(d\) is part of the data.

\begin{maintheorem}[Morita invariance; Subsection~\ref{ss:proof-A}]\label{main:morita}
Fix an integer \(d\). Let \(\TT\) and \(\UU\) be proper dg-categories admitting Serre functors. Then a \(\kk\)-linear derived Morita equivalence \(\TT\simeq\UU\) induces an equivalence between the diagrams defining \(\eta^e_n(\TT;P_\TT)\) and \(\eta^e_n(\UU;P_\UU)\) for all \(n,e\geq1\), including the Frobenius and restriction maps. In particular, we have
\[
 \hgt^e_d(\TT)=\hgt^e_d(\UU)
 \qquad(e\geq1).
\]
Moreover, if \(\pi_{-d}(\eta^e_n)\) is injective, then so is \(\pi_{-d}(\eta^e_{n+1})\).
\end{maintheorem}

Theorems~\ref{main:morita} and~\ref{main:global} rely on the functoriality of polygonic THH with respect to laced functors due to Harpaz--Nikolaus--Saunier (Theorem~\ref{thm:laced-functoriality}). A construction is given in the thesis of Saunier \cite[Theorem 6.30]{Sau25} and in his lecture notes \cite[Corollary 7.3.14]{Sau26}, and complete proofs are announced in \cite{HNS24,HNRS26}; see Remark~\ref{rem:dependence}. Theorem~\ref{main:local} below does not depend on it.

\subsection{Comparison with the classical invariants}
For a scheme \(X\) over \(\Fp\), we denote by \(\hgt(X)\) the quasi-\(F\)-split height of \(X\) in the sense of Yobuko (Definition~\ref{def:yobuko}). For a Noetherian \(F\)-finite Gorenstein local \(\Fp\)-algebra \(R\) and \(e\geq1\), we denote by \(\hgt^e(R)\) the smallest integer \(n\) such that \(R\) is \(n\)-quasi-\(F^e\)-split in the sense of \cite{TWY24}, in the form given in \cite[Definition 2.1]{Yos25} (Definition~\ref{def:yoshikawa}).

\begin{maintheorem}[Proper Gorenstein schemes; Subsection~\ref{ss:proof-B} and Corollary~\ref{cor:derived-geometric}]\label{main:global}
Let \(X\) be a nonempty proper Gorenstein scheme over \(\kk\) of pure dimension \(d\), and put \(P_X:=S_{\Perf(X)}[-d]\). Then, for every \(n\geq1\), the scheme \(X\) is \(n\)-quasi-\(F\)-split if and only if the map
\[
 \pi_{-d}\bigl(\eta^1_n(\Perf(X);P_X)\bigr)
\]
is injective. In particular, we have
\[
 \hgt^1_d(\Perf(X))=\hgt(X).
\]
Moreover, if \(Y\) is another such scheme and \(\Perf(X)\simeq\Perf(Y)\) as \(\kk\)-linear dg-categories, then \(\dim X=\dim Y\) and \(\hgt(X)=\hgt(Y)\).
\end{maintheorem}

Theorem~\ref{main:global} applies to singular Gorenstein schemes and to smooth proper varieties with nontrivial canonical bundle. The Serre functor cannot be replaced by the identity coefficient: for an ordinary smooth projective curve \(X\) of genus at least two, the Frobenius acts bijectively on \(H^1(X,\OO_X)\), so that the construction with the identity coefficient gives an injective map, whereas \(X\) is not quasi-\(F\)-split (Remark~\ref{rem:role-Serre}).

\begin{maintheorem}[Complete local Gorenstein rings; Subsection~\ref{ss:proof-C}]\label{main:local}
Let \((R,\hm,\kk)\) be a nonzero complete Noetherian Gorenstein local \(\kk\)-algebra with residue field \(\kk\), and put \(d=\dim R\). Let \(\TT_R:=\Perf_{\hm}(R)\) be the dg-category of perfect \(R\)-complexes supported at \(\hm\). Then \(\TT_R\) is proper, and its Serre functor is \([d]\); in particular, \(P_{\TT_R}\simeq\id\). Moreover, for all \(n,e\geq1\), the ring \(R\) is \(n\)-quasi-\(F^e\)-split if and only if the map
\[
 \pi_{-d}\bigl(\eta^e_n(\TT_R;\id)\bigr)
\]
is injective. In particular, we have
\[
 \hgt^e_d(\Perf_{\hm}(R))=\hgt^e(R).
\]
\end{maintheorem}

The ring \(R\) in Theorem~\ref{main:local} is automatically \(F\)-finite, since it is a quotient of \(\kk[[x_1,\ldots,x_N]]\) by the Cohen structure theorem. Theorem~\ref{main:local} applies to complete hypersurfaces \(R=\kk[[x_1,\ldots,x_N]]/(f)\), where \(f\) is a nonzero element of the maximal ideal. In this case, we have \(\hgt^e(R)=\hgt^e_{N-1}(\Perf_\hm(R))\), and the Fedder-type criterion of Yoshikawa \cite[Theorem A]{Yos25} computes this categorical invariant. For example, if \(p=2\), \(\kk\) is algebraically closed and \(R=\kk[[x,y,z]]/(z^2+x^2y+xy^n)\) with \(n\geq2\), then \(\hgt^1_2(\Perf_\hm(R))=\lceil\log_2n\rceil\) by Theorem~\ref{main:local} and \cite[Theorem 4.31]{KTY25}; see Example~\ref{ex:RDP}. The dg-category \(\TT_R\) consists of perfect complexes with finite length cohomology; it is not the singularity category.

\subsection{Relation to previous work}
Yobuko \cite[Theorem 4.5]{Yob19} shows that the quasi-\(F\)-split height of a Calabi--Yau variety coincides with the height of its Artin--Mazur formal group. Antieau--Bragg \cite[Proposition 1.1 and Theorem 5.1]{AB22} prove that this formal group, and hence its height, is a derived invariant. Therefore the derived invariance of the quasi-\(F\)-split height of Calabi--Yau varieties, including K3 surfaces, is already known.

For smooth projective varieties, the consequence \(\hgt(X)=\hgt(Y)\) of Theorem~\ref{main:global} is new only in a limited range. If \(\omega_X\) is ample or anti-ample, then \(X\cong Y\) by Bondal--Orlov \cite{BO01}. If the Kodaira dimension of \(X\) is positive, then so is that of \(Y\) (see \cite[Section 6.1]{Huy06}), and neither \(X\) nor \(Y\) is quasi-\(F\)-split by \cite[Proposition 3.14]{KTTWYY25}. Thus Theorem~\ref{main:global} gives new information for smooth projective varieties of nonpositive Kodaira dimension which are neither Calabi--Yau nor determined by their derived categories, and for singular Gorenstein schemes. In a related direction, Krah--Vial \cite{KV25} prove that global \(F\)-regularity is a derived invariant for normal Gorenstein projective varieties.

Our computation relies on the theory of Witt vectors with coefficients due to Dotto--Krause--Nikolaus--Patchkoria \cite{DKNP25} and on the polygonic spectra of \cite{KMN23}. The localization of THH with coefficients, which we need for schemes, is provided by the theory of laced categories due to Harpaz--Nikolaus--Saunier \cite{HNS24}, and for THH and TR without coefficients we use the localization theorems of Blumberg--Mandell \cite{BM12}. Liu--Liu \cite[Theorems 3.9 and 3.11]{LL26} prove that coefficient Morita equivalences induce equivalences of THH with coefficients and of the associated genuine restriction systems; Theorem~\ref{thm:laced-functoriality}, due to Harpaz--Nikolaus--Saunier, provides the corresponding functoriality for polygonic spectra. On the classical side, Lemma~\ref{lem:global-duality} and Proposition~\ref{prop:local-dual-test} are versions of the cohomological criteria in \cite[Lemmas 2.10 and 3.13]{KTTWYY25}, and dualizing complexes on Witt rings are studied in \cite[Section 9]{KTTWYY24}.

The global comparison (Theorem~\ref{main:global}) concerns only the case \(e=1\), whereas the local comparison (Theorem~\ref{main:local}) holds for all \(e\geq1\) and uses Matlis duality over \(W_n(R)\). If \(e>1\), then the classical pushout \(Q^e_{X,n}\) is only a \(W_n\OO_X\)-module, and a global comparison would require duality for \(W_n\OO_X\)-modules. This leads to the following question.

\begin{question}\label{q:global-e}
Does Theorem~\ref{main:global} extend to all \(e\geq1\), with respect to the definition of quasi-\(F^e\)-splitting for schemes in \cite{TWY24}?
\end{question}

We do not discuss quasi-\(F\)-regularity in this paper. Its definition \cite{TWY24} involves a multiplier \(c\) which runs over all non-zero-divisors, and Theorem~\ref{main:morita} does not provide the corresponding structure on dg-categories.

\subsection{Organization of the paper}
In Section~\ref{sec:classical}, we fix conventions and recall the definitions of quasi-\(F\)-splitting and quasi-\(F^e\)-splitting. In Section~\ref{sec:tr}, we recall finite TR with coefficients and its functoriality for laced categories, and we deduce its invariance under derived Morita equivalence (Corollary~\ref{cor:structured-morita}). In Section~\ref{sec:height}, we define the invariants \(\hgt^e_d\) and prove Theorem~\ref{main:morita}. In Section~\ref{sec:line-Witt}, we compute the zeroth homotopy groups of TR with coefficients in a line bundle. In Section~\ref{sec:descent}, we prove Zariski descent for TR with coefficients in a line bundle, and we compute the lowest homotopy group of derived global sections. Using these results, we prove Theorem~\ref{main:global} in Section~\ref{sec:global} and Theorem~\ref{main:local} in Section~\ref{sec:local}, where we also give examples.

\subsection*{Use of AI}
The proof strategy of this paper is due to the author. In working out the details of the arguments, the author benefited greatly from ChatGPT Astra 6. Claude Opus 5.5 was of great help in writing this paper, in particular with its English and grammar. The author takes full responsibility for all the contents of this paper.

\section{Conventions and classical definitions}\label{sec:classical}

\subsection{Conventions}
Throughout this paper, we fix a prime \(p\) and a perfect field \(\kk\) of characteristic \(p\). All tensor products of bimodules are derived, and we often omit \(\LL\) from the notation. We regard a dg-algebra as an associative ring spectrum via the Eilenberg--MacLane functor. We always consider THH relative to the sphere spectrum.

For a complex \(K\), we use the convention \(\pi_i(K)=H^{-i}(K)\), so that the shift \([1]\) corresponds to the suspension. A spectrum is connective if its negative homotopy groups vanish. For a finite group \(G\), we denote by \((-)^{hG}\), \((-)_{hG}\) and \((-)^{tG}\) the homotopy fixed points, the homotopy orbits and the Tate construction, respectively.

For an exact endofunctor \(P\) of a stable \(\infty\)-category \(\CC\), we consider the associated bimodule
\[
 M_P(x,y):=\Map_{\CC}(x,P(y)).
\]
With this convention, the bimodule associated with \(-\otimes L\) for a line bundle \(L\) is \(L\). If \(\CC\simeq\Perf(A)\), we use the same symbol for an endofunctor and the corresponding \(A\)-bimodule. Following \cite{HNS24}, we call a pair \((\CC,M)\) of a stable \(\infty\)-category \(\CC\) and a \(\CC\)-bimodule \(M\) a laced category.

\subsection{Witt vectors and quasi-\texorpdfstring{\(F\)}{F}-splitting}
Let \(A\) be an \(\Fp\)-algebra. We denote by \(W_n(A)\) the ring of \(p\)-typical Witt vectors of length \(n\), by \([a]\in W_n(A)\) the Teichm\"uller lift of \(a\in A\), and by
\[
 \res\colon W_n(A)\longrightarrow W_{n-1}(A),\qquad V\colon W_{n-1}(A)\longrightarrow W_n(A)
\]
the restriction and the Verschiebung. Let \(F_A\colon A\to A\) be the absolute Frobenius. The induced ring endomorphism \(W_n(F_A)\) of \(W_n(A)\) is given by \((a_0,\ldots,a_{n-1})\mapsto(a_0^p,\ldots,a_{n-1}^p)\) and preserves the length. We distinguish it from the Frobenius \(F\colon W_{n+1}(A)\to W_n(A)\), which decreases the length and corresponds to the Frobenius of finite TR (Proposition~\ref{prop:connective-Witt}). For a \(W_n(A)\)-module \(M\) and \(e\geq1\), we denote by \(F^e_*M\) the \(W_n(A)\)-module obtained by restriction of scalars along \(W_n(F_A)^e\).

\begin{definition}\label{def:classical-Q}
Let \(n,e\geq1\). We define a \(W_n(A)\)-module \(Q^e_{A,n}\) and a \(W_n(A)\)-linear map \(\Phi^e_{A,n}\colon A\to Q^e_{A,n}\) by the pushout diagram of \(W_n(A)\)-modules
\[
\begin{CD}
W_n(A) @>{W_n(F_A)^e}>> F_*^eW_n(A)\\
@V{\res^{n-1}}VV @VVV\\
A @>{\Phi^e_{A,n}}>> Q^e_{A,n},
\end{CD}
\]
where \(A\) is regarded as a \(W_n(A)\)-module via \(\res^{n-1}\). We put \(Q_{A,n}:=Q^1_{A,n}\) and \(\Phi_{A,n}:=\Phi^1_{A,n}\).
\end{definition}

If \(e\geq2\), then \(Q^e_{A,n}\) is not an \(A\)-module in general; for \(e=1\), we have the following description.

\begin{lemma}\label{lem:classical-quotient}
The \(W_n(A)\)-module structure on \(Q_{A,n}\) factors through \(\res^{n-1}\colon W_n(A)\to A\), and there is a natural isomorphism
\[
 Q_{A,n}\simeq F_*\bigl(W_n(A)/pW_n(A)\bigr).
\]
Under this isomorphism, the \(A\)-module structure on \(Q_{A,n}\) and the map \(\Phi_{A,n}\) are given by
\[
 a\cdot\overline w=\overline{[a]^p w},
 \qquad
 \Phi_{A,n}(a)=\overline{[a]^p}.
\]
\end{lemma}
\begin{proof}
Since \(\res^{n-1}\colon W_n(A)\to A\) is surjective with kernel \(VW_{n-1}(A)\), the pushout is the quotient of \(F_*W_n(A)\) by the image of \(VW_{n-1}(A)\) under \(W_n(F_A)\). Since \(A\) has characteristic \(p\), we have
\[
 W_n(F_A)\bigl(VW_{n-1}(A)\bigr)=pW_n(A).
\]
Hence the elements of \(VW_{n-1}(A)\) act on the quotient by zero, and the module structure factors through \(A\). The formulas follow by evaluating on Teichm\"uller lifts.
\end{proof}

This construction is compatible with localization, and hence it sheafifies on \(\Fp\)-schemes. For an \(\Fp\)-scheme \(X\), we denote by \(Q_{X,n}\) and \(\Phi_{X,n}\colon\OO_X\to Q_{X,n}\) the resulting \(\OO_X\)-module and map.

\begin{definition}\label{def:yobuko}
Let \(X\) be an \(\Fp\)-scheme and \(n\geq1\). We say that \(X\) is \(n\)-quasi-\(F\)-split if \(\Phi_{X,n}\) admits an \(\OO_X\)-linear retraction. The quasi-\(F\)-split height of \(X\) is defined by
\[
 \hgt(X):=\inf\{n\geq1\mid X\text{ is }n\text{-quasi-}F\text{-split}\}.
\]
This definition agrees with \cite[Definition 4.1]{Yob19}.
\end{definition}

\begin{definition}\label{def:yoshikawa}
Let \(R\) be a Noetherian \(F\)-finite Gorenstein local \(\Fp\)-algebra, and let \(W_n\omega_R\) be the canonical (dualizing) module of \(W_n(R)\) (see \cite[Section 9]{KTTWYY24}). Following \cite[Sections 3.1 and 3.5]{TWY24} in the form of \cite[Definition 2.1]{Yos25}, we say that \(R\) is \(n\)-quasi-\(F^e\)-split if the map
\begin{equation}\label{eq:yoshikawa-Hom}
(\Phi^e_{R,n})^*\colon
\Hom_{W_n(R)}(Q^e_{R,n},W_n\omega_R)
\longrightarrow
\Hom_{W_n(R)}(R,W_n\omega_R)
\end{equation}
is surjective. We denote by \(\hgt^e(R)\) the smallest such integer \(n\), where we put \(\hgt^e(R)=\infty\) if no such \(n\) exists.
\end{definition}

\section{Finite TR with coefficients}\label{sec:tr}

\subsection{Polygonic spectra and finite TR}
For a ring spectrum \(A\) and an \(A\)-bimodule \(P\), we put
\[
 \XX(A;P)_r:=\THH(A;P^{\otimes_A r})\qquad(r\geq1).
\]
The cyclic bar construction gives a \(C_r\)-action on this spectrum, and the Tate diagonal gives natural maps
\[
 \varphi_{p,r}\colon\XX(A;P)_r\longrightarrow
 \XX(A;P)_{pr}^{tC_p},
\]
which are equivariant with respect to the residual \(C_r\)-actions. By \cite[Theorem 6.31]{KMN23}, these data form a polygonic spectrum \(\XX(A;P)\), functorially in the \(A\)-bimodule \(P\). We only use the part indexed by powers of \(p\).

\begin{definition}\label{def:finite-TR}
For \(m\geq1\), we put \(T_m:=\{1,p,\ldots,p^{m-1}\}\) and define
\[
 \TR^m(A;P):=\TR_{T_m}(\XX(A;P)),
\]
where the right-hand side is the finite TR of polygonic spectra with respect to the truncation set \(T_m\), which is defined as a mapping spectrum in \cite[Section 3]{KMN23}.
\end{definition}

Explicitly, writing \(\XX_r=\XX(A;P)_r\), the spectrum \(\TR^m(A;P)\) is the equalizer of two maps
\[
 \prod_{0\leq i<m}\XX_{p^i}^{hC_{p^i}}
 \rightrightarrows
 \prod_{1\leq i<m}\bigl(\XX_{p^i}^{tC_p}\bigr)^{hC_{p^{i-1}}},
\]
one of which is induced by the canonical maps and the other by the maps \(\varphi_p\). In particular, finite TR preserves finite limits of polygonic spectra.

The inclusion \(T_m\subset T_{m+1}\) induces the restriction map \(\res\colon\TR^{m+1}(A;P)\to\TR^m(A;P)\). The shift of \(\XX(A;P)\) by \(p\) is identified with \(\XX(A;P^{\otimes_Ap})\), and the Frobenius of finite TR gives a map
\begin{equation}\label{eq:F-coeff}
 F\colon\TR^{m+1}(A;P)\longrightarrow
 \TR^m(A;P^{\otimes_A p}).
\end{equation}
These maps commute with the restriction maps and can be iterated.

\subsection{Functoriality and Morita invariance}
We use the following functoriality of polygonic THH for laced categories, due to Harpaz--Nikolaus--Saunier.

\begin{theorem}[Harpaz--Nikolaus--Saunier]\label{thm:laced-functoriality}
The following assertions hold.
\begin{enumerate}
\item The construction \((\CC,M)\mapsto\bigl(\THH(\CC;M^{\circ r})\bigr)_{r\geq1}\), together with the cyclic actions and the Tate diagonals, extends to a functor \(\XX\) from the \(\infty\)-category of laced categories to the \(\infty\)-category of polygonic spectra in the sense of \cite{KMN23}; see \cite[Theorem 6.30]{Sau25} and \cite[Corollary 7.3.14]{Sau26}. For \(M=\id\), it recovers the cyclotomic structure on \(\THH(\CC)\).
\item If \(A\) is a connective ring spectrum and \(P\) is a connective \(A\)-bimodule, then \(\XX(\Perf(A),P)\) is naturally equivalent to the polygonic spectrum \(\XX(A;P)\) of \cite[Theorem 6.31]{KMN23}; see \cite[Section 7.2]{Sau25}.
\item The structure maps satisfy \(\varphi_{p,r}(\CC,M)=\varphi_{p,1}(\CC,M^{\circ r})\), naturally in \((\CC,M)\); see \cite[Section 6.4]{Sau25} and \cite[Section 7.2]{Sau26}. In particular, the shift of \(\XX(\CC,M)\) by \(p\) is identified with \(\XX(\CC,M^{\circ p})\), naturally in \((\CC,M)\).
\end{enumerate}
\end{theorem}

\begin{remark}\label{rem:dependence}
The construction in \cite[Section 6.4]{Sau25} is presented as a sketch, statement (ii) is stated without proof in \cite[Section 7.2]{Sau25} (see also \cite[Proposition 7.16]{Sau25}), and \cite{Sau26} consists of lecture notes. A genuine refinement of (i) is stated in \cite[Theorem 7.18]{Sau25}, and complete proofs are announced in \cite[Introduction]{HNS24} and \cite{HNRS26}. Theorems~\ref{main:morita} and~\ref{main:global} depend on Theorem~\ref{thm:laced-functoriality}. Theorem~\ref{main:local} does not: its coefficient is equivalent to the identity, and if \(\TR^m(\TT_R;\id)\) is defined via the cyclotomic spectrum \(\THH(\TT_R)\), then its proof only uses the functoriality of \(\XX(A;P)\) in \(P\) for a fixed dg-algebra \(A\) \cite[Theorem 6.31]{KMN23}, the comparison of TR of cyclotomic spectra with the classical TR \cite{KMN23}, and the invariance and localization theorems for THH and TR \cite{BM12}.
\end{remark}

For a laced category \((\CC,P)\), we write \(\THH(\CC;P)\) for the first level of \(\XX(\CC,P)\) and \(\TR^m(\CC;P):=\TR_{T_m}(\XX(\CC,P))\). By Theorem~\ref{thm:laced-functoriality}, the restriction and Frobenius maps of finite TR are natural with respect to laced functors.

\begin{corollary}\label{cor:structured-morita}
Let \(A\) and \(B\) be dg-algebras, and let \(M\) be an \(A\)-\(B\)-bimodule and \(N\) a \(B\)-\(A\)-bimodule such that
\[
 M\otimes_BN\simeq A,\qquad N\otimes_AM\simeq B
\]
as bimodules. For an \(A\)-bimodule \(P\), we put \(P':=N\otimes_AP\otimes_AM\). Then there is an equivalence of polygonic spectra
\[
 \XX(\Perf(A),P)\simeq\XX(\Perf(B),P'),
\]
which induces equivalences between the finite TR spectra compatible with the Frobenius and restriction maps.
\end{corollary}
\begin{proof}
The functor \(-\otimes_AM\colon\Perf(A)\to\Perf(B)\) is an equivalence, and we have a natural equivalence \((-\otimes_AP)\otimes_AM\simeq(-\otimes_AM)\otimes_BP'\) since \(M\otimes_BN\simeq A\). Hence it underlies an equivalence of laced categories \((\Perf(A),P)\simeq(\Perf(B),P')\). By Theorem~\ref{thm:laced-functoriality}(i), it induces an equivalence of polygonic spectra, which is compatible with the identifications of the shifts by \(p\) by Theorem~\ref{thm:laced-functoriality}(iii). Since finite TR and its restriction and Frobenius maps are natural in polygonic spectra, we obtain the claim.
\end{proof}

In particular, for a dg-category \(\TT\) and an exact endofunctor \(P\) of \(\TT\), the spectra \(\THH(\TT;P)\) and \(\TR^m(\TT;P)\), together with their restriction and Frobenius maps, are invariant under derived Morita equivalence. If \(\TT\simeq\Perf(A)\) with \(A\) and \(P\) connective, then they are computed by \(\XX(A;P)\) by Theorem~\ref{thm:laced-functoriality}(ii).

\subsection{Connectivity and the zeroth homotopy groups}
\begin{proposition}\label{prop:connective-Witt}
Let \(A\) be a connective ring spectrum and \(P\) a connective \(A\)-bimodule. Then \(\TR^m(A;P)\) is connective, and there is a natural isomorphism
\[
 \pi_0\TR^m(A;P)
 \simeq W_m(\pi_0A;\pi_0P)
\]
compatible with restriction, Frobenius and Verschiebung. Here \(W_m(\pi_0A;\pi_0P)\) denotes the Witt vectors of length \(m\) with coefficients in \(\pi_0P\) in the sense of \cite{DKNP25}.
\end{proposition}
\begin{proof}
The fiber of \(\res\colon\TR^{m+1}(A;P)\to\TR^m(A;P)\) is given by the homotopy orbits of the corresponding cyclic tensor product. It is connective, since all tensor factors are connective and taking homotopy orbits preserves connectivity. Hence the first assertion follows by induction on \(m\); compare \cite[Lemma 3.14 and the proof of Proposition 5.2]{KMN23}.

The second assertion is \cite[Theorem 2.5]{DKNP25}, and the compatibility with the operators follows from \cite[Proposition 2.13]{DKNP25} together with the compatibility with restriction in its construction. Our index \(m\) is the number of elements of the truncation set; thus \(W_m\) in this paper is denoted by \(W_{\langle p^{m-1}\rangle}\) in \cite{DKNP25}.
\end{proof}

\section{The invariants and their Morita invariance}\label{sec:height}

\subsection{Definition}
Let \(\TT\) be a dg-category and \(P\) an exact endofunctor of \(\TT\). We define \(\cQ^e_n(\TT;P)\) by \eqref{eq:intro-Q}, and we denote by
\[
 \eta^e_n(\TT;P)\colon\THH(\TT;P)\longrightarrow\cQ^e_n(\TT;P)
\]
the composite of the inclusion of the second summand and the canonical map to the cofiber. Thus we have a pushout square of spectra
\[
\begin{CD}
\TR^{n+e}(\TT;P) @>{F^e}>> \TR^n(\TT;P^{\otimes p^e})\\
@V{\res^{n+e-1}}VV @VVV\\
\THH(\TT;P) @>{\eta^e_n(\TT;P)}>> \cQ^e_n(\TT;P).
\end{CD}
\]

\begin{definition}\label{def:height}
Let \(d\) be an integer and \(e\geq1\). We define
\[
 \hgt^e_d(\TT;P)
 :=\inf\{n\geq1\mid\pi_{-d}\bigl(\eta^e_n(\TT;P)\bigr)\text{ is injective}\}.
\]
If \(\TT\) is proper and admits a Serre functor \(S_\TT\), we put
\[
 P_{\TT}:=S_{\TT}[-d],
 \qquad
 \hgt^e_d(\TT):=\hgt^e_d(\TT;P_{\TT}).
\]
\end{definition}

\begin{lemma}\label{lem:monotone}
There are maps
\[
 \cQ^e_{n+1}(\TT;P)\longrightarrow\cQ^e_n(\TT;P)
\]
induced by the restriction maps, whose composite with \(\eta^e_{n+1}(\TT;P)\) is \(\eta^e_n(\TT;P)\). In particular, if \(\pi_{-d}(\eta^e_n(\TT;P))\) is injective, then \(\pi_{-d}(\eta^e_{n'}(\TT;P))\) is injective for every \(n'\geq n\).
\end{lemma}
\begin{proof}
We apply the restriction maps to the two TR terms and the identity to the THH term. Since \(F^e\circ\res=\res\circ F^e\) and the two composites of restriction maps to THH agree, we obtain a map of pushout diagrams, and hence the required map. If \(\pi_{-d}(\eta^e_n)\) is injective, then any element killed by \(\pi_{-d}(\eta^e_{n+1})\) is also killed by the composite \(\pi_{-d}(\eta^e_n)\), and hence it is zero.
\end{proof}

\begin{proposition}\label{prop:coefficient-height-morita}
A derived Morita equivalence \(\TT\simeq\UU\) which carries \(P\) to \(P'\) induces an equivalence of arrows
\[
 \eta^e_n(\TT;P)\simeq\eta^e_n(\UU;P').
\]
In particular, we have \(\hgt^e_d(\TT;P)=\hgt^e_d(\UU;P')\) for all \(e\) and \(d\).
\end{proposition}
\begin{proof}
By Corollary~\ref{cor:structured-morita}, both arrows in the diagram defining \(\cQ^e_n\) are identified. Since pushouts are functorial, the pushouts and the maps from the THH terms are also identified. Hence the two maps on \(\pi_{-d}\) are conjugate by isomorphisms, and one of them is injective if and only if the other is.
\end{proof}

\subsection{\texorpdfstring{Proof of Theorem~\ref{main:morita}}{Proof of Theorem A}}\label{ss:proof-A}
Let \(E\colon\TT\simeq\UU\) be the given equivalence, and choose dg-algebras \(A\) and \(B\) with \(\TT\simeq\Perf(A)\) and \(\UU\simeq\Perf(B)\). Then \(E\) is given by an invertible bimodule, and it carries the \(A\)-bimodule \(A^\vee\) to the \(B\)-bimodule \(B^\vee\), since both are the \(\kk\)-linear duals of the corresponding mapping bimodules. Hence \(E\circ S_{\TT}\simeq S_{\UU}\circ E\) and \(E\circ P_{\TT}\simeq P_{\UU}\circ E\), and the first assertion follows from Proposition~\ref{prop:coefficient-height-morita}. The last assertion is Lemma~\ref{lem:monotone}.
\qed

\begin{remark}
The definition of \(\hgt^e_d\) does not use a distinguished unit element in \(\HH_0\) or in \(\TR\). Such an element is not invariant under Morita equivalence: under the Morita equivalence \(\HH_0(M_r(\kk))\simeq\HH_0(\kk)=\kk\), the class of the identity matrix corresponds to \(r\).
\end{remark}

\section{Witt vectors with coefficients in a line bundle}\label{sec:line-Witt}

Let \(X\) be a Noetherian \(\Fp\)-scheme and \(L\) an invertible \(\OO_X\)-module. We denote by
\[
 W_m(\OO_X;L)
\]
the sheaf of Witt vectors of length \(m\) with coefficients in \(L\), regarded as a symmetric bimodule, which is obtained by gluing the local constructions of \cite{DKNP25}.

\begin{lemma}\label{lem:Witt-line}
The sheaf \(W_m(\OO_X;L)\) is an invertible \(W_m\OO_X\)-module. More precisely, let \(\ell\) and \(\ell'=g\ell\) be local generators of \(L\), where \(g\) is a unit. Then the transition map from the trivialization given by \(\ell'\) to the trivialization given by \(\ell\) is
\[
 (a_0,\ldots,a_{m-1})
 \longmapsto(g a_0,g^p a_1,\ldots,g^{p^{m-1}}a_{m-1}),
\]
that is, the multiplication by \([g]\).
\end{lemma}
\begin{proof}
If \(L\) is trivialized by \(\ell\), then the tensor powers appearing in the construction of Witt vectors with coefficients are identified with the corresponding powers of \(\ell\), and hence \(W_m(\OO_X;L)\) is identified with \(W_m\OO_X\). The multiplication by \(g\) on \(L\) acts on the \(i\)-th coordinate, which lies in \(L^{\otimes p^i}\), by \(g^{p^i}\); this is the multiplication by the Teichm\"uller lift \([g]\). Since \([gh]=[g][h]\), these transition maps satisfy the cocycle condition, and the local modules glue.
\end{proof}

\begin{remark}\label{rem:Witt-divisorial}
By Lemma~\ref{lem:Witt-line}, \(W_m(\OO_X;L)\) is the Teichm\"uller lift of \(L\). In particular, if \(X\) is normal and \(L=\OO_X(D)\) for a Cartier divisor \(D\), then \(W_m(\OO_X;L)\) is isomorphic to the Witt divisorial sheaf \(W_m\OO_X(D)\) of \cite[Subsection 2.6]{KTTWYY25}, which is locally generated by \([f]^{-1}\) if \(D=\operatorname{div}(f)\) locally. For \(L=\omega_X\), this sheaf should be distinguished from the dualizing sheaf \(W_m\omega_X\) of \(W_m\OO_X\).
\end{remark}

We have the maps
\[
 \res^n\colon W_{n+1}(\OO_X;L)\longrightarrow L,
 \qquad
 F\colon W_{n+1}(\OO_X;L)\longrightarrow W_n(\OO_X;L^{\otimes p}).
\]

\begin{proposition}\label{prop:line-pushout}
The pushout of the above two maps in the category of sheaves of abelian groups is naturally isomorphic to
\[
 Q_{X,n}\otimes_{\OO_X}L,
\]
where \(Q_{X,n}\) is equipped with the \(\OO_X\)-module structure of Lemma~\ref{lem:classical-quotient}. Under this isomorphism, the induced map from \(L\) is \(\Phi_{X,n}\otimes L\).
\end{proposition}
\begin{proof}
The map \(\res^n\) is surjective with kernel
\[
 VW_n(\OO_X;L^{\otimes p}).
\]
For general coefficients, the composite \(FV\) is the norm map of the cyclic action on tensor powers. Since this action is trivial on tensor powers of an invertible module regarded as a symmetric bimodule, we have \(FV=p\) on \(W_n(\OO_X;L^{\otimes p})\). Hence the pushout is the sheaf of abelian groups
\[
 W_n(\OO_X;L^{\otimes p})/p.
\]
We equip this quotient with the \(\OO_X\)-module structure given by \(a\cdot w=[a]^p w\). If \(\ell\) is a local generator of \(L\), then this quotient is identified with \(Q_{X,n}\) by Lemma~\ref{lem:classical-quotient}. If \(\ell'=g\ell\) is another local generator, then the transition map is the multiplication by \([g]^p\) by Lemma~\ref{lem:Witt-line}, which is exactly the transition map of \(Q_{X,n}\otimes L\) with respect to the \(\OO_X\)-module structure of Lemma~\ref{lem:classical-quotient}.

Moreover, with respect to the generator \(\ell\), an element \(a\ell\) of \(L\) is sent to the class of \([a]^p\ell^{\otimes p}\), which corresponds to \(\Phi_{X,n}(a)\otimes\ell\). Hence these identifications glue, and the induced map is \(\Phi_{X,n}\otimes L\).
\end{proof}

\begin{corollary}\label{cor:pi0-line-Q}
Let \(U=\Spec A\) be an affine open subscheme of \(X\), and regard \(L(U)\) as a symmetric \(A\)-bimodule. Then the spectra \(\THH(A;L(U))\) and \(\cQ^1_n(A;L(U))\) are connective, and there are natural isomorphisms
\[
 \pi_0\THH(A;L(U))\simeq L(U),\qquad
 \pi_0\cQ^1_n(A;L(U))\simeq(Q_{X,n}\otimes L)(U),
\]
compatible with the restriction to smaller affine open subschemes. Under these isomorphisms, the map \(\pi_0(\eta^1_n)\) corresponds to \(\Phi_{X,n}\otimes L\).
\end{corollary}
\begin{proof}
By Proposition~\ref{prop:connective-Witt}, all the terms in the diagram defining \(\cQ^1_n(A;L(U))\) are connective, and their zeroth homotopy groups are given by Witt vectors with coefficients. Since a cofiber of a map of connective spectra is connective and its \(\pi_0\) is the cokernel of the induced map on \(\pi_0\), the assertion follows from Proposition~\ref{prop:line-pushout}.
\end{proof}

\section{Zariski descent and the lowest homotopy group}\label{sec:descent}

\subsection{Zariski descent}
Let \(X\) be a quasi-compact separated scheme and \(L\) a line bundle on \(X\). For an open subscheme \(U\subset X\), we denote by \(P_U\) the endofunctor \(-\otimes L|_U\) of \(\Perf(U)\).

\begin{proposition}\label{prop:line-descent}
The polygonic spectra \(\XX(\Perf(-),P_{(-)})\), the spectra \(\TR^m(\Perf(-);P_{(-)})\) and the spectra \(\cQ^e_n(\Perf(-);P_{(-)})\) satisfy Zariski descent; that is, they send a cover of a quasi-compact open subscheme by two quasi-compact open subschemes to a Mayer--Vietoris square. Moreover, if \(X\) is Noetherian of finite Krull dimension, then these spectra for \(X\) are identified with the derived global sections of the associated Zariski sheaves of spectra.
\end{proposition}
\begin{proof}
For quasi-compact open subschemes \(V\subset U\), the restriction functor \(\Perf(U)\to\Perf(V)\) commutes with \(-\otimes L\), and hence it underlies a laced functor \((\Perf(U),P_U)\to(\Perf(V),P_V)\). By Theorem~\ref{thm:laced-functoriality}, we obtain a functor from quasi-compact open subschemes of \(X\) to polygonic spectra. Since the Tate construction is exact, finite limits of polygonic spectra are computed levelwise, and hence it suffices to show that \(\THH(\Perf(-);P^r)\) satisfies Mayer--Vietoris for every \(r\geq1\).

Let \(j\colon U\hookrightarrow X\) be a quasi-compact open immersion, and put \(Z=X\setminus U\). For every \(r\geq1\), the endofunctor \(P_X^r=-\otimes L^{\otimes r}\) preserves \(\Perf_Z(X)\), and the endofunctor induced on \(\Perf(U)\) is \(P_U^r\); on Ind-categories, this follows from
\[
 j^*\bigl(j_*E\otimes L^{\otimes r}\bigr)
 \simeq E\otimes j^*L^{\otimes r}.
\]
Thus the localization sequence \(\Perf_Z(X)\to\Perf(X)\to\Perf(U)\) gives, for each \(r\), a laced Verdier sequence up to idempotent completion. By \cite[Proposition 5.18 and Corollary 5.19]{HNS24}, THH with these coefficients sends it to a fiber sequence.
If \(X=U\cup V\), then the dg-categories of perfect complexes supported on \(X\setminus U\) in \(X\) and in \(V\) are equivalent, compatibly with the coefficients, by excision. Comparing the localization sequences, we obtain a Mayer--Vietoris square for \(\THH(\Perf(-);P^r)\).

Since finite TR preserves finite limits of polygonic spectra, and since the pushout defining \(\cQ^e_n\) is exact and its Frobenius term is natural by Theorem~\ref{thm:laced-functoriality}(iii), the spectra \(\TR^m\) and \(\cQ^e_n\) also satisfy Mayer--Vietoris. Finally, if \(X\) is Noetherian of finite Krull dimension, a finite affine open cover and its intersections compute the associated sheaf by iterated Mayer--Vietoris squares, and the convergence of the descent spectral sequence (Lemma~\ref{lem:bottom}) identifies the constructions for \(X\) with the derived global sections.
\end{proof}

\subsection{The lowest homotopy group}
\begin{lemma}\label{lem:bottom}
Let \(X\) be a Noetherian scheme of dimension \(d\), and let \(\mathscr E\) be a sheaf of spectra on the Zariski site of \(X\) whose homotopy sheaves \(\pi_t\mathscr E\) vanish for \(t<0\). Then we have
\[
 \pi_i\R\Gamma(X,\mathscr E)=0\quad(i<-d),
 \qquad
 \pi_{-d}\R\Gamma(X,\mathscr E)
 \simeq H^d(X,\pi_0\mathscr E),
\]
and the isomorphism is natural in \(\mathscr E\).
\end{lemma}
\begin{proof}
Since the underlying topological space of \(X\) is Noetherian of Krull dimension \(d\), the \(\infty\)-topos of sheaves on \(X\) is locally of homotopy dimension at most \(d\) by \cite[Proposition 7.2.4.7 and Corollary 7.2.4.17]{HTT}, and hence it is hypercomplete by \cite[Corollary 7.2.1.12]{HTT}. Moreover, the Zariski cohomological dimension of \(X\) is at most \(d\) by \cite[\stacks{02UZ}]{Stacks}. Hence the descent spectral sequence
\[
 E_2^{s,t}=H^s(X,\pi_t\mathscr E)\Longrightarrow
 \pi_{t-s}\R\Gamma(X,\mathscr E)
\]
converges. Since \(E_2^{s,t}=0\) unless \(0\leq s\leq d\) and \(t\geq0\), there are no nonzero terms of total degree \(t-s<-d\), and the only term of total degree \(-d\) is \(E_2^{d,0}\). No nonzero differential enters this term, since its source would have \(t<0\), and no nonzero differential leaves it, since its target would have \(s>d\). Hence we obtain the claim.
\end{proof}

\begin{lemma}\label{lem:support-cube}
Let \(R\) be a Noetherian ring, let \(t_1,\ldots,t_d\in R\), and put \(I=(t_1,\ldots,t_d)\) and \(t_J=\prod_{i\in J}t_i\) for a subset \(J\subset\{1,\ldots,d\}\). Let \(J\mapsto E_J\) be a cubical diagram of connective spectra indexed by the subsets of \(\{1,\ldots,d\}\), and let \(E^{\mathrm{tot}}\) be its total fiber. Then \(\pi_iE^{\mathrm{tot}}=0\) for \(i<-d\). Moreover, if the diagram \(J\mapsto\pi_0E_J\) is isomorphic to the diagram \(J\mapsto M_{t_J}\) of localizations of an \(R\)-module \(M\), then there is a natural isomorphism
\[
 \pi_{-d}E^{\mathrm{tot}}\simeq H^d_I(M).
\]
\end{lemma}
\begin{proof}
Filtering the total fiber by the cardinality of \(J\), we obtain a spectral sequence
\[
 E_1^{s,t}=\bigoplus_{|J|=s}\pi_tE_J\Longrightarrow\pi_{t-s}E^{\mathrm{tot}},
\]
where \(t\geq0\) and \(0\leq s\leq d\). Hence \(\pi_iE^{\mathrm{tot}}=0\) for \(i<-d\), and \(\pi_{-d}E^{\mathrm{tot}}\) is the \(d\)-th cohomology of the complex \(E_1^{\bullet,0}\). Under the assumption, this complex is the \v{C}ech complex of \(M\) with respect to \(t_1,\ldots,t_d\), whose \(d\)-th cohomology is \(H^d_I(M)\).
\end{proof}

\section{Proper Gorenstein schemes}\label{sec:global}

Let \(X\) be a nonempty proper Gorenstein scheme over \(\kk\) of pure dimension \(d\). Then \(\omega_X\) is invertible, and \(\omega_X[d]\) is a dualizing complex of \(X\). Since \(X\) is proper, \(\Perf(X)\) is proper. By \cite[\stacks{09IS}]{Stacks}, there is a perfect complex which generates \(D_{\mathit{QCoh}}(\OO_X)\); since the compact objects of \(D_{\mathit{QCoh}}(\OO_X)\) are the perfect complexes \cite[\stacks{09M1}]{Stacks}, it is a classical generator of \(\Perf(X)\). The Serre functor of \(\Perf(X)\) is
\[
 S_X=-\otimes\omega_X[d];
\]
see \cite[\stacks{0FY8}]{Stacks}. Hence the coefficient in Definition~\ref{def:height} is \(P_X=-\otimes\omega_X\).

\subsection{\texorpdfstring{Proof of Theorem~\ref{main:global}}{Proof of Theorem B}}\label{ss:proof-B}
\begin{proposition}\label{prop:global-arrow}
The map \(\pi_{-d}(\eta^1_n(\Perf(X);P_X))\) is naturally identified with the map
\begin{equation}\label{eq:global-top-map}
 H^d(X,\omega_X)
 \xrightarrow{H^d(\Phi_{X,n}\otimes\omega_X)}
 H^d(X,Q_{X,n}\otimes\omega_X).
\end{equation}
\end{proposition}
\begin{proof}
By Proposition~\ref{prop:line-descent}, the source and the target of \(\eta^1_n(\Perf(X);P_X)\), together with the map between them, are identified with the derived global sections of the corresponding Zariski sheaves of spectra. By Theorem~\ref{thm:laced-functoriality}(ii) on affine open subschemes and Corollary~\ref{cor:pi0-line-Q}, the homotopy sheaves of these sheaves vanish in negative degrees, and their zeroth homotopy sheaves and the map between them are identified with
\[
 \omega_X
 \xrightarrow{\Phi_{X,n}\otimes\omega_X}
 Q_{X,n}\otimes\omega_X.
\]
Hence the assertion follows from Lemma~\ref{lem:bottom} and its naturality.
\end{proof}

\begin{lemma}\label{lem:global-duality}
Let \(Q\) be a coherent \(\OO_X\)-module and \(\Phi\colon\OO_X\to Q\) an \(\OO_X\)-linear map. Then the induced map
\[
 H^d(X,\omega_X)\longrightarrow H^d(X,Q\otimes\omega_X)
\]
is injective if and only if \(\Phi\) admits an \(\OO_X\)-linear retraction.
\end{lemma}
\begin{proof}
Since \(X\) is an equidimensional proper Cohen--Macaulay scheme, Serre duality gives a natural isomorphism
\[
 H^d(X,\mathcal F)^\vee\simeq\Hom_X(\mathcal F,\omega_X)
\]
for every coherent \(\OO_X\)-module \(\mathcal F\); see \cite[\stacks{0FVZ}]{Stacks}. Applying it to \(\mathcal F=Q\otimes\omega_X\) and \(\mathcal F=\omega_X\), and using the invertibility of \(\omega_X\), we identify the dual of the above map with
\[
 \Phi^*\colon\Hom_X(Q,\OO_X)\longrightarrow\Hom_X(\OO_X,\OO_X).
\]
Since all the groups are finite-dimensional over \(\kk\), the above map is injective if and only if \(\Phi^*\) is surjective. If \(\Phi^*\) is surjective, then a preimage of \(\id_{\OO_X}\) is a retraction of \(\Phi\). Conversely, if \(\psi\) is a retraction of \(\Phi\), then \(\Phi^*(u\circ\psi)=u\) for every endomorphism \(u\) of \(\OO_X\).
\end{proof}

\begin{proof}[Proof of the first assertion of Theorem~\ref{main:global}]
Since \(X\) is of finite type over a perfect field, the Frobenius of \(X\) is finite, and hence \(Q_{X,n}\) is a coherent \(\OO_X\)-module. By Proposition~\ref{prop:global-arrow} and Lemma~\ref{lem:global-duality}, the map \(\pi_{-d}(\eta^1_n(\Perf(X);P_X))\) is injective if and only if \(\Phi_{X,n}\) admits an \(\OO_X\)-linear retraction, that is, \(X\) is \(n\)-quasi-\(F\)-split. Since this holds for every \(n\), we obtain \(\hgt^1_d(\Perf(X))=\hgt(X)\).
\end{proof}

\subsection{Dimension and derived invariance}
Let \(\TT\) be a nonzero proper dg-category with a Serre functor \(S\) and a classical generator \(G\). Following Elagin--Lunts \cite[Definition 5.3]{EL21}, we put
\[
 e_-(G,G'):=\inf\{i\mid\Hom_\TT(G,G'[i])\neq0\},\qquad
 e_+(G,G'):=\sup\{i\mid\Hom_\TT(G,G'[i])\neq0\},
\]
and define the upper and lower Serre dimensions by
\[
 \overline{\Sdim}(\TT):=\limsup_{m\to\infty}\frac{-e_-(G,S^mG)}{m},
 \qquad
 \underline{\Sdim}(\TT):=\liminf_{m\to\infty}\frac{-e_+(G,S^mG)}{m}.
\]
They do not depend on the choice of \(G\) \cite{EL21}, and they are invariant under \(\kk\)-linear equivalences, which commute with Serre functors.

\begin{proposition}\label{prop:dimension}
Let \(X\) be as above. Then \(\overline{\Sdim}(\Perf(X))=\underline{\Sdim}(\Perf(X))=d\). In particular, \(d\) is invariant under \(\kk\)-linear equivalences of these dg-categories.
\end{proposition}
\begin{proof}
For smooth projective \(X\), this is \cite[Lemma 5.6]{EL21}; the same argument applies. Let \(G\) be a classical generator of \(\Perf(X)\). Since \(G^\vee\otimes G\) is perfect, its cohomology sheaves are concentrated in a bounded range of degrees, and this range does not change after tensoring with \(\omega_X^{\otimes m}\). Since coherent cohomology on \(X\) vanishes in degrees greater than \(d\), there are integers \(a\leq b\), independent of \(m\), such that
\[
 H^j\bigl(X,G^\vee\otimes G\otimes\omega_X^{\otimes m}\bigr)=0
 \quad\text{unless }a\leq j\leq b.
\]
Since
\[
 \Hom_X(G,S_X^mG[i])
 \simeq
 H^{md+i}\bigl(X,G^\vee\otimes G\otimes\omega_X^{\otimes m}\bigr),
\]
this group vanishes unless \(a-md\leq i\leq b-md\). On the other hand, at least one of these groups is nonzero, since otherwise we would have \(S_X^mG=0\). Hence \(a-md\leq e_-(G,S_X^mG)\leq e_+(G,S_X^mG)\leq b-md\), and dividing by \(m\), we obtain the claim.
\end{proof}

\begin{corollary}\label{cor:derived-geometric}
Let \(X\) and \(Y\) be nonempty proper Gorenstein schemes over \(\kk\) of pure dimension. If there is a \(\kk\)-linear derived Morita equivalence \(\Perf(X)\simeq\Perf(Y)\), then
\[
 \dim X=\dim Y,\qquad \hgt(X)=\hgt(Y).
\]
In particular, quasi-\(F\)-splitting and \(F\)-splitting are preserved by such equivalences.
\end{corollary}
\begin{proof}
The assertion follows from Proposition~\ref{prop:dimension}, Theorem~\ref{main:morita} and the first assertion of Theorem~\ref{main:global}, since \(X\) is \(F\)-split if and only if \(\hgt(X)=1\).
\end{proof}

\subsection{The role of the Serre functor}
\begin{remark}\label{rem:role-Serre}
Assume that \(\kk\) is algebraically closed, and let \(X\) be an ordinary smooth projective curve over \(\kk\) of genus at least two. Since the Kodaira dimension of \(X\) is one, \(X\) is not quasi-\(F\)-split by \cite[Proposition 3.14]{KTTWYY25}, and hence \(\hgt^1_1(\Perf(X))=\infty\) by Theorem~\ref{main:global}. On the other hand, the Frobenius acts bijectively on \(H^1(X,\OO_X)\). The argument of Proposition~\ref{prop:global-arrow} with \(L=\OO_X\) identifies \(\pi_{-1}(\eta^1_1(\Perf(X);\id))\) with the Frobenius on \(H^1(X,\OO_X)\), which is injective. Hence the identity coefficient does not recover the quasi-\(F\)-split height, and the shift of the Serre functor is necessary in Theorem~\ref{main:global}.
\end{remark}

\section{Complete local Gorenstein rings}\label{sec:local}

Throughout this section, let \((R,\hm,\kk)\) be a nonzero complete Noetherian Gorenstein local \(\kk\)-algebra with residue field \(\kk\), and put \(d=\dim R\). As mentioned in the introduction, \(R\) is \(F\)-finite.

\subsection{The dg-category of supported perfect complexes}
\begin{proposition}\label{prop:local-Serre}
The dg-category \(\TT_R=\Perf_{\hm}(R)\) is proper, has a classical generator, and admits the Serre functor \([d]\).
\end{proposition}
\begin{proof}
For \(K,L\in\TT_R\), the complex \(\RHom_R(K,L)\) is a perfect \(R\)-complex with finite length cohomology, and hence \(\TT_R\) is proper. A Koszul complex \(K_0\) on a system of parameters of \(R\) is a classical generator of \(\TT_R\); see \cite[\stacks{0A9A}, \stacks{0A9B}]{Stacks}. We put \(A:=\RHom_R(K_0,K_0)\), so that \(\TT_R\simeq\Perf(A)\).

Since \(R\) is Gorenstein, \(R[d]\) is a normalized dualizing complex of \(R\). For an \(R\)-module \(M\) of finite length, we have \(\operatorname{Ext}^i_R(M,R)=0\) for \(i\neq d\) and \(\operatorname{Ext}^d_R(M,R)\simeq\Hom_R(M,E_R(\kk))\) by local duality \cite[\stacks{0A81}]{Stacks}. Since \(\kk\subset R\) is a coefficient field, we have \(E_R(\kk)\simeq\colim_q\Hom_{\kk}(R/\hm^q,\kk)\), and the evaluation at \(1\) identifies \(\Hom_R(M,E_R(\kk))\) with \(\Hom_\kk(M,\kk)\); see \cite[\stacks{0A7Q}]{Stacks}. Hence we obtain a natural equivalence
\[
 \RHom_\kk(C,\kk)\simeq\RHom_R(C,R[d])
\]
of functors on \(R\)-complexes \(C\) with finite length cohomology. Since the \(A\)-bimodule structure on \(A\) is \(R\)-linear, applying this to \(C=A\) gives an equivalence of \(A\)-bimodules \(A^\vee\simeq\RHom_R(A,R)[d]\). Since \(K_0\) is perfect, the trace pairing \(A\otimes_RA\to R\), \((u,v)\mapsto\tr(uv)\), induces an equivalence of \(A\)-bimodules \(A\simeq\RHom_R(A,R)\). Hence \(A^\vee\simeq A[d]\) as \(A\)-bimodules, that is, the Serre functor of \(\TT_R\) is \([d]\).
\end{proof}

Consequently, we have \(P_{\TT_R}=S_{\TT_R}[-d]\simeq\id\).

\subsection{The zeroth homotopy groups}
For a discrete commutative \(\Fp\)-algebra \(A\), we put \(\cQ^e_n(A):=\cQ^e_n(\Perf(A);\id)\).

\begin{proposition}\label{prop:local-pi0}
The spectrum \(\cQ^e_n(A)\) is connective, and there is a natural isomorphism
\[
 \pi_0\cQ^e_n(A)\simeq Q^e_{A,n}
\]
under which the map \(\pi_0(\eta^e_n)\) corresponds to \(\Phi^e_{A,n}\).
\end{proposition}
\begin{proof}
By Proposition~\ref{prop:connective-Witt}, the three spectra \(\TR^{n+e}(A)\), \(\TR^n(A)\) and \(\THH(A)\) in the diagram defining \(\cQ^e_n(A)\) are connective, and their zeroth homotopy groups are \(W_{n+e}(A)\), \(W_n(A)\) and \(A\), respectively.

Since \(A\) has characteristic \(p\), the Frobenius \(F^e\colon W_{n+e}(A)\to W_n(A)\) is given by
\[
 (a_0,\ldots,a_{n+e-1})
 \longmapsto(a_0^{p^e},\ldots,a_{n-1}^{p^e}).
\]
Hence it factors as the surjection \(\res^e\colon W_{n+e}(A)\to W_n(A)\) followed by \(W_n(F_A)^e\). The other map \(\res^{n+e-1}\) factors as \(\res^{n-1}\circ\res^e\). In particular, both maps vanish on \(\ker(\res^e)\).

Since all the spectra are connective, \(\pi_0\) of the pushout is the pushout of abelian groups of the zeroth homotopy groups. Dividing the source by the common kernel \(\ker(\res^e)\), we identify it with the underlying abelian group of the pushout of \(W_n(A)\)-modules in Definition~\ref{def:classical-Q}, compatibly with the maps from \(A\). This also equips \(\pi_0\cQ^e_n(A)\) with the \(W_n(A)\)-module structure of \(Q^e_{A,n}\).
\end{proof}

\subsection{Local cohomology}
\begin{proposition}\label{prop:local-arrow}
The map \(\pi_{-d}(\eta^e_n(\TT_R;\id))\) is naturally identified with the map
\begin{equation}\label{eq:local-top-map}
 H^d_{\hm}(R)
 \xrightarrow{H^d_{\hm}(\Phi^e_{R,n})}
 H^d_{\hm}(Q^e_{R,n}).
\end{equation}
\end{proposition}
\begin{proof}
Since the coefficient is the identity, the polygonic spectra here come from cyclotomic spectra, and finite TR agrees with the classical TR \cite{KMN23}. Put \(U=\Spec R\setminus\{\hm\}\). By the localization theorem \cite[Theorems 1.1 and 1.3]{BM12}, \(\THH(\TT_R)\) is identified with the fiber of \(\THH(R)\to\THH(U)\), compatibly with the cyclotomic structures. Since finite TR and the cofiber defining \(\cQ^e_n\) are exact, \(\cQ^e_n(\TT_R)\) is also identified with the fiber of \(\cQ^e_n(R)\to\cQ^e_n(U)\).

Choose a system of parameters \(t_1,\ldots,t_d\) of \(R\). Then the principal open subschemes \(D(t_i)\) cover \(U\). By Zariski descent \cite{BM12}, \(\THH(\TT_R)\) and \(\cQ^e_n(\TT_R)\) are identified with the total fibers of the cubes \(J\mapsto\THH(R_{t_J})\) and \(J\mapsto\cQ^e_n(R_{t_J})\), which consist of connective spectra. For THH, the diagram of zeroth homotopy groups is the diagram \(J\mapsto R_{t_J}\) of localizations of \(R\).

For \(\cQ^e_n\), by Proposition~\ref{prop:local-pi0}, the diagram of zeroth homotopy groups is \(J\mapsto Q^e_{R_{t_J},n}\), which is the diagram of localizations of the \(W_n(R)\)-module \(Q^e_{R,n}\) at the Teichm\"uller lifts \([t_J]\). Indeed, we have \(W_n(R)_{[t]}\simeq W_n(R_t)\), and localizing \(F^e_*W_n(R)\) at \([t]\) is the same as localizing it at \([t]^{p^e}\). Since the kernel of \(W_n(R)\to R\) is nilpotent, the ideal \(([t_1],\ldots,[t_d])\) of \(W_n(R)\) has the maximal ideal as its radical.

Applying Lemma~\ref{lem:support-cube} to the Noetherian ring \(W_n(R)\) (Lemma~\ref{lem:Witt-CM}) and the elements \([t_1],\ldots,[t_d]\), we obtain
\[
 \pi_{-d}\THH(\TT_R)\simeq H^d_{\hm}(R),
 \qquad
 \pi_{-d}\cQ^e_n(\TT_R)\simeq H^d_{\hm}(Q^e_{R,n}),
\]
where both local cohomology groups are taken over \(W_n(R)\) with respect to its maximal ideal; the first one agrees with the usual local cohomology of \(R\). By naturality, we obtain the claim.
\end{proof}

\subsection{\texorpdfstring{Duality over \(W_n(R)\) and the proof of Theorem~\ref{main:local}}{Duality over W\_n(R) and the proof of Theorem C}}\label{ss:proof-C}
For the convenience of the reader, we include a proof of the following ring-theoretic facts.

\begin{lemma}\label{lem:Witt-CM}
The ring \(W_n(R)\) is a complete Noetherian Cohen--Macaulay local ring of dimension \(d\). Moreover, \(Q^e_{R,n}\) is a finite \(W_n(R)\)-module.
\end{lemma}
\begin{proof}
Put \(W_0(R)=0\). The Verschiebung filtration of \(W_n(R)\) has graded pieces
\[
 V^iW_{n-i}(R)/V^{i+1}W_{n-i-1}(R)\simeq F_*^iR\qquad(0\leq i\leq n-1).
\]
On the \(i\)-th graded piece, \([a]\) acts as the multiplication by \(a^{p^i}\); in particular, the action of \(W_n(R)\) on each graded piece factors through the surjection \(W_n(R)\to R\).

Firstly, we prove that \(W_n(R)\) is Noetherian by induction on \(n\), using the exact sequence
\[
 0\longrightarrow V^{n-1}R\longrightarrow W_n(R)
 \longrightarrow W_{n-1}(R)\longrightarrow0.
\]
The kernel is isomorphic to \(F_*^{n-1}R\), which is a finite \(R\)-module since \(R\) is \(F\)-finite. Hence it is a Noetherian \(W_n(R)\)-module. Since \(W_{n-1}(R)\) is a Noetherian ring by the induction hypothesis, it is also a Noetherian \(W_n(R)\)-module, and hence \(W_n(R)\) is Noetherian.

Secondly, since the kernel of \(W_n(R)\to R\) is nilpotent, \(W_n(R)\) is local, its spectrum is homeomorphic to \(\Spec R\), and its dimension is \(d\). On each graded piece \(F^i_*R\), the topology defined by the maximal ideal of \(W_n(R)\) coincides with the \(\hm\)-adic topology, since the ideal generated by the \(p^i\)-th powers of the elements of \(\hm\) is \(\hm\)-primary. These finite \(R\)-modules are complete. Since completion is exact for finite modules over a Noetherian ring, the above exact sequence and induction on \(n\) show that \(W_n(R)\) is complete.

Thirdly, let \(t_1,\ldots,t_d\) be a system of parameters of \(R\), which is a regular sequence since \(R\) is Cohen--Macaulay. On the \(i\)-th graded piece, \([t_1],\ldots,[t_d]\) act as \(t_1^{p^i},\ldots,t_d^{p^i}\), which form a regular sequence on \(F^i_*R\). Since a sequence which is regular on modules \(M'\) and \(M''\) is regular on every extension of \(M''\) by \(M'\), we see by induction on the filtration that \([t_1],\ldots,[t_d]\) is a regular sequence on \(W_n(R)\). Hence \(W_n(R)\) is Cohen--Macaulay.

Finally, \(F_*^eW_n(R)\) has a finite filtration whose graded pieces are \(F_*^{e+i}R\), which are finite \(W_n(R)\)-modules. Hence \(F_*^eW_n(R)\) is a finite \(W_n(R)\)-module. Since \(R\) is also a finite \(W_n(R)\)-module, so is the pushout \(Q^e_{R,n}\).
\end{proof}

\begin{proposition}\label{prop:local-dual-test}
The map \eqref{eq:local-top-map} is injective if and only if \(R\) is \(n\)-quasi-\(F^e\)-split.
\end{proposition}
\begin{proof}
Let \(\omega_{W_n(R)}=W_n\omega_R\) be a canonical module of \(W_n(R)\), and let
\[
 D_n(-):=\Hom_{W_n(R)}(-,E_{W_n(R)}(\kk))
\]
be the Matlis duality over \(W_n(R)\). By Lemma~\ref{lem:Witt-CM} and local duality \cite[\stacks{0A81}]{Stacks}, we have a natural isomorphism
\[
 D_n H^d_{\hm}(M)\simeq\Hom_{W_n(R)}(M,\omega_{W_n(R)})
\]
for every finite \(W_n(R)\)-module \(M\). The local cohomology modules \(H^d_\hm(M)\) are Artinian, and Matlis duality is exact and faithful on Artinian modules. Hence, applying \(D_n\) to \eqref{eq:local-top-map}, we see that \eqref{eq:local-top-map} is injective if and only if the map
\[
 \Hom_{W_n(R)}(Q^e_{R,n},W_n\omega_R)
 \longrightarrow
 \Hom_{W_n(R)}(R,W_n\omega_R)
\]
is surjective. This is the condition in Definition~\ref{def:yoshikawa}.
\end{proof}

\begin{proof}[Proof of Theorem~\ref{main:local}]
By Proposition~\ref{prop:local-Serre}, the dg-category \(\TT_R\) is proper and its Serre functor is \([d]\), and hence \(P_{\TT_R}\simeq\id\). By Proposition~\ref{prop:local-arrow}, the map \(\pi_{-d}(\eta^e_n(\TT_R;\id))\) is identified with the map \eqref{eq:local-top-map} of local cohomology groups. By Proposition~\ref{prop:local-dual-test}, this map is injective if and only if \(R\) is \(n\)-quasi-\(F^e\)-split. Since this holds for every \(n\), we obtain \(\hgt^e_d(\Perf_\hm(R))=\hgt^e(R)\).
\end{proof}

For \(e=1\), Proposition~\ref{prop:local-dual-test} and \cite[Lemma 3.13]{KTTWYY25} show that \(\hgt^1(R)\) coincides with the quasi-\(F\)-split height \(\hgt(R)\) of \(\Spec R\) in the sense of Definition~\ref{def:yobuko}.

\begin{remark}[A dg-algebra model of a singularity]
For a Koszul complex \(K_0\) as in the proof of Proposition~\ref{prop:local-Serre}, the dg-algebra \(A_R=\RHom_R(K_0,K_0)\) is proper, and we have \(\Perf(A_R)\simeq\Perf_{\hm}(R)\). Hence Theorem~\ref{main:local} expresses the classical heights in terms of a proper dg-algebra model of \(R\), and the result does not depend on the choice of \(K_0\).
\end{remark}

\subsection{Examples}
\begin{example}\label{ex:dual-numbers}
Let \(A=\kk[\epsilon]/(\epsilon^2)\). Then \(A\) is a zero-dimensional complete Gorenstein local ring, and \(\Perf(A)=\Perf_{(\epsilon)}(A)\) is a proper dg-category which is not smooth and whose Serre functor is the identity. For all \(n,e\geq1\), we have \(\Phi^e_{A,n}(\epsilon)=\overline{[\epsilon^{p^e}]}=0\) by the definition of the pushout. Since \(d=0\), the map \eqref{eq:local-top-map} is \(\Phi^e_{A,n}\) itself, and hence \(\hgt^e_0(\Perf(A))=\infty\) for all \(e\geq1\) by Theorem~\ref{main:local}. The same argument applies to every zero-dimensional Gorenstein local \(\kk\)-algebra \(A\) which is not a field: a nonzero element \(s\) of the socle satisfies \(s^2=0\), and hence \(\Phi^e_{A,n}(s)=0\). Thus finite heights greater than one only appear in positive dimension.
\end{example}

\begin{example}\label{ex:RDP}
Assume that \(\kk\) is algebraically closed of characteristic \(2\). For \(n\geq2\), let
\[
 R=\kk[[x,y,z]]/(z^2+x^2y+xy^n),
\]
which is a rational double point. By \cite[Theorem 4.31]{KTY25}, the quasi-\(F\)-split height of \(R\) is \(\lceil\log_2n\rceil\). Hence, by Theorem~\ref{main:local} and the comparison for \(e=1\) above, we have
\[
 \hgt^1_2(\Perf_\hm(R))=\lceil\log_2n\rceil.
\]
In particular, the invariant \(\hgt^1_2\) takes every positive integer as a value on proper dg-categories of the form \(\Perf_\hm(R)\).
\end{example}

\bibliographystyle{amsalpha}
\bibliography{references}
\end{document}